\documentclass[11pt,twoside]{article}
\usepackage{a4wide}
\usepackage{amssymb}
\usepackage{amsmath}
\usepackage{enumerate}
\newtheorem{theorem}{Theorem}
\newtheorem{mtp}{Theorem (Mass Transference Principle)}

\newtheorem{corollary}{Corollary}
\newtheorem{lemma}{Lemma}

\renewcommand{\Bbb}[1]{\mathbb{#1}}
\newcommand{\N}{{\Bbb N}}         % natural numbers
\newcommand{\R}{{\Bbb R}}         % real numbers
\newcommand{\Z}{{\Bbb Z}}         % integer numbers

\newcommand{\cH}{{\cal H}}

\newcommand{\ie}{{\it i.e.}\/ }

\renewcommand{\le}{\leqslant}
\renewcommand{\ge}{\geqslant}

\begin{document}

\title{On a problem of K. Mahler: Diophantine approximation \\ and  Cantor sets  }

\author{Jason Levesley \\ {\small\sc York } \and  Cem Salp  \\ {\small\sc York}  \and  Sanju Velani\footnote{Royal
Society University Research Fellow} \\ {\small\sc York}}

\maketitle

\date

\centerline{{\it Dedicated to Maurice Dodson on his retirement --
finally! }}

\begin{abstract}  Let $K$ denote the middle third Cantor set and
${\cal A}:= \{ 3^n : n = 0,1,2, \ldots \} $.  Given a
 real, positive function $\psi$ let  $ W_{\cal
A}(\psi)$ denote the set of real numbers $x$ in the unit interval
for which there exist infinitely many $(p,q) \in \Z \times {\cal
A} $ such that $ |x - p/q| < \psi(q)  $.  The analogue of  the
Hausdorff measure version of the Duffin-Schaeffer conjecture is
established for  $ W_{\cal A}(\psi) \cap K  $. One of the
consequences of this is that  there exist very well approximable
numbers, other than Liouville numbers, in $K$ -- an assertion
attributed to K. Mahler.
\end{abstract}

\bigskip


%\newpage

\section{Introduction}

\subsection{A problem of K. Mahler \label{prob}}

A real number $x$ is said to be {\em very well approximable}  if
there exists some $\epsilon > 0$ such that
\begin{equation} \label{vwa}
|x - p/q|  \ < \  q^{-(2+ \epsilon)} \hspace{6mm} {\rm \ for \
infinitely \ many \ } \hspace{6mm} (p,q) \in \Z \times \N \ .
\end{equation}
 Note that in view of
Dirichlet's theorem, if  $ \epsilon = 0 $ then  every real number
satisfies the above inequality. The set  ${\cal W}$ of very well
approximable numbers is a set of Lebesgue measure zero but
nevertheless is large in the sense that it has maximal Hausdorff
dimension; i.e. $\dim {\cal W}= 1 $.  A real number $x$ is said to
be a Liouville number if (\ref{vwa}) is satisfied for all $\epsilon
> 0 $. It is well known that the set ${\cal L}$ of Liouville numbers is uncountable
and of zero Hausdorff dimension; i.e. $\dim {\cal L} = 0 $.

\medskip

Throughout, let $K$ denote the standard middle third Cantor set.
Thus, $K$ is precisely the set of  real numbers in the unit interval
whose base three expansions are free of the digit one. It is well
known that $K$ is a set of Lebesgue measure zero and
$$
\dim K  \ =  \ \gamma \ := \frac{\log 2}{\log 3 } \ .
$$

\vspace{2ex}

The following assertion is attributed  to K. Mahler  -- see Problem
35 in \cite[\S10.2]{bug} and {\rm Remark (iii)} below. %references to
%Mahler within.

\vspace{2ex}

\noindent {\bf Mahler's Assertion. \ }  {\em There exists very
well approximable numbers, other than Liouville numbers, in the
middle third Cantor set; i.e. $$ ({\cal W}\setminus {\cal L })\cap
K  \ \neq  \ \emptyset \ . $$}

\vspace{2ex}

It is rather surprising that this claim  remains unproved to this
day. A simple consequence of our main result -- Theorem \ref{t1}
below, is that
\begin{equation} \label{walk}
\dim ( \, ({\cal W}\setminus {\cal L })\cap K  )  \ \geq \ \gamma/ 2
\ .
\end{equation}

\noindent Clearly, this `strongly' implies the assertion of Mahler.
In fact, we shall prove  a lot more than (\ref{walk}). We show that
there exist real numbers in $K$ with any prescribed exact order --
see \S\ref{eo}.

%At the end of the paper
In \S\ref{explicit} we give explicit examples of irrational
numbers satisfying Mahlers assertion.

\vspace*{2ex}

\noindent{\em {Remarks:}}
\begin{enumerate}[\em (i)]
\item It is easy to see that the real
number, $$ \xi :=  2 \, \sum_{n=1}^{\infty}   \, 3^{-n!}
$$ is a Liouville number and therefore very well approximable.
Furthermore, $\xi $ clearly lies in $K$ since only the digits $0$
and $2$ appear in its base $3$ expansion. Hence, ${\cal W}\cap K
\neq   \emptyset $ and it is therefore  natural to exclude Liouville
numbers in the above assertion.
\item %On a slightly different note, a
A real number $x$ is said to be badly approximable if there exists a
constant $c(x)> 0 $ such that $|x -p/q| \geq c(x)/q^2 $ for all
rational $p/q$. Let ${\bf Bad}$ denote the set of badly approximable
numbers. It is well known that ${\bf Bad} \cap K  $ is non-empty
and moreover that $\dim ({\bf Bad} \cap K) = \gamma $ -- see for
example \cite{KW,KTV}.
\item We have not been able to find a direct source to Mahler in which  the
above form of the assertion is stated.  In view of {\rm Remarks
(i)} \& {\rm (ii)}, Mahler's assertion as stated above is in all
likelihood a more precise reformulation of the following problem
posed by Mahler in \cite[\S2]{mah}:  {\em How close can irrational
elements of Cantor's set be approximated by rational numbers?} In
any case, the results obtained in \S\ref{eo}  provide a
satisfactory and precise solution to this rather vague problem. In
short, the `vague' answer is: {\em However close one wishes!}
\end{enumerate}

\section{The set $W_{\cal A}(\psi)$ and our approach }

Throughout, $\psi:\R^+\to\R^+$ will denote a real, positive
function and, unless explicitly stated otherwise,   ${\cal A}:= \{
3^n : n = 0,1,2, \ldots \} $. Also
$$
W_{\cal A}(\psi):= \{x \in [0,1] :|x - p/q| < \psi(q) \ {\rm for\
i.m.\ } \; (p,q) \in \Z \times {\cal A} \} \  , $$ where `i.m'
stands for `infinitely many'.  Thus, $W_{\cal A}(\psi)$ is simply
a hybrid of the classical set $W(\psi)$ of $ \psi$--well
approximable numbers in which the denominator of the rational
approximates are restricted to the set ${\cal A }$;  in other
words if we were to put  ${\cal A} = \N $ then $ W_{\cal A}(\psi)
= W(\psi)$. In the case $\psi : r \to r^{-\tau}$ with $\tau > 0$,
we write $W_{\cal A}(\tau)$ for $W_{\cal A}(\psi)$. Thus, when
${\cal A} = \N $ then $ W_{\cal A}(\tau) = W(\tau)$ is the
classical set of $\tau$--well approximable numbers.

\subsection{The general metric theory for $W_{\cal A}(\psi) \cap K $ \label{gmt}}
Let $f$ be a dimension function and let $\cH^f$ denote the
Hausdorff $f$--measure -- see \S\ref{HM}.  In short, our aim is to
provide a complete metric theory for the set $W_{\cal A}(\psi)
\cap K $. The following result achieves this goal in that it
provides a  simple criteria for the `size' of the set $W_{\cal
A}(\psi) \cap K $ expressed in terms of the general measure
$\cH^f$.

\begin{theorem}\label{t1}
Let $f$ be a dimension function such that $r^{-\gamma}f(r)$ is
monotonic.  Then
$$
\cH^f(W_{\cal A } (\psi) \cap K )=\left\{
\begin{array}{cl}
0& {\rm if}  \qquad\displaystyle \sum_{n=1}^{\infty} \
f\!\left(\psi(3^n)\right)  \times \ (3^n)^{\gamma} \ < \ \infty\, \\[4ex]
\cH^f(K)  & {\rm if}  \qquad\displaystyle\sum_{n=1}^{\infty} \
f\!\left(\psi(3^n)\right)  \times \ (3^n)^{\gamma}   \ = \ \infty\,
\end{array}
\right..
$$
\end{theorem}

The convergence part of the above theorem is relatively
straightforward if not trivial -- see \S\ref{conv}. The main
substance is the divergent part. It is worth stressing that we do
not assume that the function $\psi$ is monotonic. Thus, within the
framework under consideration the above theorem establishes the
analogue of the general form of the   Duffin-Schaeffer conjecture
as formulated in \cite{mtp}. The fact that we have not imposed the
condition that $(p,q)=1 $ in the definition of $W_{\cal A }
(\psi)$ is irrelevant since $\phi(q) := \# \{1 \leq t \leq q:
(t,q)=1 \} = \frac23 \,  q $ for any $q \in {\cal A}$ and so   $$
\displaystyle \sum_{n=1}^{\infty} \ f\!\left(\psi(3^n)\right)
\times \ (3^n)^{\gamma}  \ \asymp \ \displaystyle
\sum_{n=1}^{\infty} \ f\!\left(\psi(3^n)\right) \times \
(\phi(3^n))^{\gamma} \ \ .
$$
For details regarding the original statement of Duffin and
Schaeffer see \cite{bug,Harman}.

\medskip

With $f: r  \to r^s$ ($s \geq 0$), an immediate consequence of
Theorem \ref{t1} is the following corollary.

\begin{corollary}\label{c1} For $\tau \geq 1 $,
$ \dim (W_{\cal A } (\tau) \cap K ) \ = \ \gamma/ \tau\ . $ In
particular, for $\tau > 1 $
$$\cH^{\gamma/\tau} (W_{\cal A } (\tau) \cap K ) \ = \ \infty \ . $$
\end{corollary}

Now if $\tau$ is strictly greater than two, then every point in
$W_{\cal A } (\tau)$ is by definition very well approximable. Thus,
$$
\dim ({\cal W} \cap K  )  \ \geq \ \gamma/ 2 \ . $$ This together
with the fact that the set ${\cal L }$ of Liouville numbers  is of
dimension zero implies (\ref{walk}). In turn, this implies the
assertion of Mahler.

\section{Sets of exact order in $ K $ \label{eo}} Recall, that
for  $\tau > 0 $  the set  $W(\tau)$ of $\tau$-well approximable
numbers  consists of real numbers such that $$|x - p/q| \ < \
q^{-\tau} \hspace{6mm} {\rm \ for \ infinitely \ many \ }
\hspace{6mm} (p,q) \in \Z \times \N \ .  $$ For a real number $x$,
its {\em exact order} $\tau(x) $ is defined as follows: $$ \tau(x)
:= \sup \{ \tau: x \in W(\tau) \} \; . $$ It follows from
Dirichlet's theorem that $\tau(x)\ge 2$ for all $x \in \R$. For
$\alpha \geq 2$, let $E(\alpha)$ denote the set of numbers with {\em
exact order} $\alpha$; that is $$ E(\alpha) := \{ x \in \R : \tau(x)
= \alpha \} \; . $$ Thus, $ E(\alpha)$ consists of real numbers with
`order' of rational approximation sandwiched between $ \alpha -
\epsilon $ and $\alpha + \epsilon $ where $\epsilon > 0 $ is
arbitrarily small. The set $E(\alpha)$ is equivalent to the set of
real numbers $x$ for which Mahler's function $\theta_1(x)$ is equal
to $ \alpha - 1 $; the general function $\theta_n (x) $ is central
to Mahler's classification of transcendental numbers (see
\cite{BS,bug,Gut}). A simple consequence of Theorem \ref{t1} is the
following result which trivially implies Mahler's assertion and
provides a precise solution to the problem discussed in {\rm Remark
(iii)}  of \S\ref{prob}.

\vspace*{2ex}

\begin{corollary}\label{c2} For $\alpha \geq 2 $,
$ \dim ( E(\alpha) \cap K ) \ \geq \ \gamma/ \alpha\ . $ In
particular, for $\alpha > 2 $
$$\cH^{\gamma/\alpha} ( E(\alpha) \cap K ) \ = \ \infty \ . $$
\end{corollary}

\vspace*{2ex}

\noindent{\em Proof. \ }  When $\alpha =2 $, we have  that $
E(\alpha) = \R $ and so $\dim ( E(\alpha) \cap K ) = \dim K :=
\gamma $. Thus, without loss of generality assume that $\alpha >
2$. Let
$$
\psi_1 : r \to r^{-\alpha} \hspace*{6ex}  {\rm and \ } \hspace*{6ex}
\psi_2 : r \to r^{-\alpha} (\log r )^{-\frac{2 \, \alpha}{\gamma}} \
.
$$
It is easily verified that
\begin{equation}
\label{exin} W_{\cal A } (\psi_1) \backslash W_{\cal A } (\psi_2) \
\subset \ E(\alpha)  \ .
\end{equation}
Now with   $f: r \to r^{\gamma/\alpha}$ we have that
$$
\sum_{n=1}^{\infty} \ f\!\left(\psi_1(3^n)\right) \times
(3^n)^{\gamma} \ = \ \sum_{n=1}^{\infty} 1 \ = \ \infty
$$
and
$$
\sum_{n=1}^{\infty} \ f\!\left(\psi_2(3^n)\right) \times
(3^n)^{\gamma} \ \asymp \ \sum_{n=1}^{\infty} \frac{1}{n^{2}} \ < \
\infty \ \ .
$$
By Theorem \ref{t1},
$$
\cH^{\gamma/\alpha} (W_{\cal A } (\psi_1) \cap K ) \ = \ \infty
\hspace*{4ex} {\rm and \ } \hspace*{4ex}  \cH^{\gamma/\alpha}
(W_{\cal A } (\psi_2) \cap K ) \ = \ 0 \ \ .
$$
Thus $$ \cH^{\gamma/\alpha} (W_{\cal A } (\psi_1) \backslash W_{\cal
A } (\psi_2) \cap K ) \ = \ \infty \ ,  $$ which together with
(\ref{exin}) implies the desired measure and dimension statements.
\hfill $\spadesuit$

\vspace*{2ex}

%For $\alpha \geq \frac{1}{2} (\sqrt{5} +3) $,
Explicit examples
of irrational numbers in $E(\alpha) \cap K$ are given in
\S\ref{explicit}.

 \vspace*{2ex}

\noindent{\em Remark.  } It is evident from the above proof that the
statement of Corollary \ref{c2} does not in anyway utilize the full
power of Theorem \ref{t1}. The argument outlined above can be
modified to establish much stronger `exact order' statements in the
spirit of those in \cite{exactBDV}. Essentially, the set $E(\alpha)
$ in the statement of Corollary \ref{c2} can be replaced by sets $
E(\psi,\phi)$ consisting of real numbers whose rational
approximation properties  are sandwiched between two functions
$\psi$ and $\phi$ with $\phi$ in some sense `smaller' than $\psi$.
In short,  $ E(\psi,\phi):= W_{\cal A } (\psi) \backslash W_{\cal A
} (\phi)$.

\section{Preliminaries}

\subsection{Hausdorff measures \label{HM}}

In this section we give a brief account of Hausdorff measures. For
further details see \cite{falc,MAT}. A {\em dimension function} $f
\, : \, \R^+ \to \R^+ $ is a monotonic, non-decreasing function
such that $f(r)\to 0$ as $r\to 0 \, $. The Hausdorff $f$--measure
with respect to the dimension function $f$ will be denoted
throughout by ${\cal H}^{f}$ and is defined as follows. Suppose
$F$ is  a subset  of $\R^k$. For $\rho
> 0$, a countable collection $ \left\{B_{i} \right\} $ of balls in
$\R^k$ with radius  $r(B_i) \leq \rho $ for each $i$ such that $F
\subset \bigcup_{i} B_{i} $ is called a {\em $ \rho $-cover for
$F$}.
 For a dimension
function $f$ define $$
 {\cal H}^{f}_{\rho} (F) \, = \, \inf \ \sum_{i} f(r(B_i)),
$$
where the infimum is taken over all $\rho$-covers of $F$. The {\it
Hausdorff $f$--measure} $ {\cal H}^{f} (F)$ of $F$ with respect to
the dimension function $f$ is defined by   $$ {\cal H}^{f} (F) :=
\lim_{ \rho \rightarrow 0} {\cal H}^{f}_{\rho} (F) \; = \;
\sup_{\rho > 0 } {\cal H}^{f}_{\rho} (F) \; . $$

%A simple consequence of the definition of $ {\cal H}^f $ is the
%following useful fact.
%
%\begin{lemma}
%If $ \, f$ and $g$ are two dimension functions such that the ratio
%$f(r)/g(r) \to 0 $ as $ r \to 0 $, then ${\cal H}^{f} (F) =0 $
%whenever ${\cal H}^{g} (F) < \infty $. \label{dimfunlemma}
%\end{lemma}
%

In the case that  $f(r) = r^s$ ($s \geq 0$), the measure $ \cH^f $
is the usual $s$--dimensional Hausdorff measure $\cH^s $ and the
Hausdorff dimension $\dim F$ of a set $F$ is defined by $$ \dim \,
F \, := \, \inf \left\{ s : {\cal H}^{s} (F) =0 \right\} = \sup
\left\{ s : {\cal H}^{s} (F) = \infty \right\} . $$ In particular
when $s$ is an integer, $\cH^s$ is comparable to the
$s$--dimensional Lebesgue measure. Actually, $\cH^s$ is a constant
multiple of the $s$--dimensional Lebesgue measure.
%but we shall not need this stronger statement.

\subsection{The  Mass Transference Principle \label{secmtp}}

Let $X$ be a  compact set in $\R^n$.  Suppose there exist
constants $\delta > 0 $, $0<c_1<1<c_2<\infty$ and $r_0 > 0$  such
that
\begin{equation*}\label{g}
c_1\ r^{\delta} \ \le \  \cH^{\delta}(B) \ \le \  c_2\ r^{\delta}
\ ,
\end{equation*}
for any ball $B=B(x,r)$ in $X$  with $x\in X$ and $r\le r_0$.  Next,
given a dimension function $f$ and a ball $B=B(x,r)$ we define
$$
B^f:=B(x,f(r)^{1/\delta})\,.
$$
When $f(r)=r^s$ for some $s>0$ we also adopt the notation $B^s$,
\ie $ B^s:=B^{(r\mapsto r^s)}$. Thus, by definition,
$B^{\delta}(x,r)=B(x,r)$. Given a sequence of balls $B_i$ in $X$,
$i=1,2,3,\ldots$, as usual its limsup set is
$$
\limsup_{i\to\infty}B_i:=\bigcap_{j=1}^\infty\ \bigcup_{i\ge j}B_i \
.
$$
By definition,  $\limsup_{i\to\infty}B_i$ is precisely the set of
points in $X$ which   lie in infinitely many balls $B_i$.

\vspace*{1ex}

The following  Mass Transference Principle allows us to transfer
$\cH^{\delta}$-measure theoretic statements for $\limsup$ subsets
of $X$ to general $\cH^f$-measure theoretic statements.

\begin{mtp}%[Mass Transference Principle]
%\label{thm3}
Let $X$  be as above and let $\{B_i\}_{i\in\N}$ be a sequence of
balls in $X$ with $r(B_i)\to 0$ as $i\to\infty$. Let $f$ be a
dimension function such that $r^{-\delta} f(x)$ is monotonic and
suppose that for any ball $B$ in $X$
\begin{equation*}\label{e:011b}
\cH^{\delta}\big(\/B\cap\limsup_{i\to\infty}B^f_i{}\,\big)=\cH^{\delta}(B)
\ .
\end{equation*}
Then, for any ball $B$ in $X$
\begin{equation*}\label{e:012b}
\cH^f\big(\/B\cap\limsup_{i\to\infty}B^{\delta}_i\,\big)=\cH^f(B)
\ .
\end{equation*}
\end{mtp}

\medskip

\noindent The theorem is essentially  Theorem~3 in \cite{mtp}. It
is  simplified  for the particular application we have in mind.

\subsection{Positive and full measure sets}

Let $X$ be a  compact set in $\R^n$ and  $\mu$ be a finite measure
supported on $X$. The measure $\mu$ is said to be {\em doubling} if
there exists a  constant $\lambda
> 1 $ such that for $x \in X$
\begin{equation}
\label{doub}
 \mu(B(x,2r)) \, \leq \, \lambda\,  \mu(B(x,r)) \ .
\end{equation}

In this section we state two measure theoretic results which will be
required during the course of establishing the divergent part of
Theorem \ref{t1}.

\begin{lemma}\label{lem1a}
Let $X$ be a  compact set in $\R^n$ and let $\mu$ be a finite
doubling measure on $X$ such that any open set is $\mu$
measurable. Let $E$ be a Borel subset of $X$. Assume that there
are constants $r_0,c>0$ such that for any ball $B$ with $r(B)<r_0$
and centre in $X$ we have that
$$\mu(E\cap B) \ \ge \  c\; \mu(B) \  .  $$ Then, $E$ has full
measure in $X$, i.e. $\mu (X \setminus E) = 0 $.
\end{lemma}

 For the proof see \cite[\S8]{BDV03}.

\begin{lemma}\label{lem2a}
Let $X$ be a  compact set in $\R^n$ and let $\mu$ be a finite
measure on $X$. Also, let  $E_n$ be a sequence of $\mu$-measurable
sets such that $\sum_{n=1}^\infty \mu(E_n)=\infty $. Then $$ \mu(
\limsup_{n \to \infty} E_n ) \; \geq \;  \limsup_{Q \to \infty}
\frac{ \left( \sum_{s=1}^{Q} \mu(E_s) \right)^2 }{ \sum_{s, t =
1}^{Q} \mu(E_s \cap E_t ) }  \  \  \ . $$
\end{lemma}

Lemma~\ref{lem2a} is relatively well known and is  a generalization
of the divergent part of the standard Borel--Cantelli lemma,  see
Lemma 5 in \cite{Sp}.

%\subsection{The  Cantor measure $\mu$  \label{mu}}

\section{Proof of Theorem \ref{t1}: the convergence part \label{conv}
}

%\subsection{The convergence part of Theorem \ref{t1}  \label{conv} }

We are given that  $f$ is a dimension function such that
\begin{equation} \label{consum}
\displaystyle \sum_{n=1}^{\infty} \ f\!\left(\psi(3^n)\right) \times
\ (3^n)^{\gamma} \ < \ \infty\ .
\end{equation}
In view of this, $ \psi(3^n)  \to 0 $ as $ n \to \infty  $. Thus
given any $\rho > 0$, there exists an integer $n_0(\rho)$  such
that
$$
\psi(3^n) \leq \rho \quad {\rm for \ all \ } \quad n \geq n_0(\rho)
\ \ .
$$
Furthermore and  without loss of generality we can assume that
$n_0(\rho) \to \infty $ as $ \rho \to 0 $.
%$$
%\psi(3^n)  \to 0  \quad {\rm as \ } \quad n \to \infty \ \ .
%$$
Now for $n \in \N $, let
\begin{eqnarray*}
 A_n  \ & :=  & \ \bigcup_{0
\ \leq \ p\  \leq \ 3^n } B( \textstyle{\frac{p}{3^n}}, \
\psi(3^{n}) )   \ \cap \ K \\ & & \\
\ & =  & \ \bigcup_{\substack{0 \ \leq \ p \  \leq \ 3^n : \\ ~ \\
\textstyle{\frac{p}{3^n} \in K  }} } B( \textstyle{\frac{p}{3^n}}, \
\psi(3^{n}) )   \ \cap \ K
\end{eqnarray*}
Then by definition,  $W_{\cal A } (\psi) \cap K = \limsup_{n \to
\infty} A_n $ and for each $m \in \N $ we have that
$$
W_{\cal A } (\psi) \cap K  \ \subset \ \bigcup_{n\geq m} A_n \ \ .
$$
It now follows from the definition of $ {\cal H}^{f} $ that for
any $\rho > 0 $,
\begin{eqnarray*}
{\cal H}^{f}_{\rho}  (W_{\cal A } (\psi) \cap K)  \ & \ll &
\sum_{n\geq n_0(\rho)} f\!\left(\psi(3^n)\right) \ \times  \ \#
\{0 \leq  p  \leq  3^n : \textstyle{\frac{p}{3^n}} \in K \} \\ & &
\\  \ & \ll &  \sum_{n\geq n_0(\rho)} f\!\left(\psi(3^n)\right)
\times (3^n)^{\gamma}   \ .
\end{eqnarray*}
In view of (\ref{consum})  and the fact that $n_0(\rho) \to \infty
$ as $ \rho \to 0 $, we have that $$ \sum_{n\geq n_0(\rho)}
f\!\left(\psi(3^n)\right) \times (3^n)^{\gamma}  \ \  \ \to \ \ \
0 \quad {\rm as \ } \ \ \rho \to 0 \ .
$$
Thus, ${\cal H}^{f} (W_{\cal A } (\psi) \cap K) = 0 $ as required.
\hfill $\spadesuit$

\section{Proof of Theorem \ref{t1}: the divergent part \label{?}
}

The divergent part of Theorem \ref{t1} constitutes the main
substance of theorem.  The proof will be split into various key
and natural  steps.

\subsection{A reduction
to the measure $\mu$ \label{div} }
%\subsection{The divergence  part of Theorem \ref{t1}: a reduction
%to the measure $\mu$ \label{div} }
%  to ${\cal H}^{\gamma} $ \label{div} }

We begin by imposing a co-primeness condition. Let $W^*_{\cal
A}(\psi) $ denote the set of $x \in [0,1] $ for which there exist
infinitely many co-prime $(p,q) \in \Z \times {\cal A}$ such that
$$  |x - p/q| < \psi(q) \ \ . $$
Trivially, $W^*_{\cal A}(\psi)  \subset W_{\cal A}(\psi) $ and so
\begin{equation} \label{imp}
 {\cal H}^{f} (W^*_{\cal A } (\psi) \cap K) = {\cal H}^{f}(K)
\ \ \Longrightarrow  \ \ {\cal H}^{f} (W_{\cal A } (\psi) \cap K) =
{\cal H}^{f} (K)  \ .
\end{equation}

Recall, that $\gamma := \dim K $. Let $\mu$ denote the restriction
of the  $\gamma$--dimensional Hausdorff measure  to $K$; that is
$$ \mu \ := \ {\cal H}^{\gamma}|_{\substack{ ~ \\ K }}  \ .
$$
It is well known that $\mu$ is a finite measure supported on $K$
and moreover  there exist constants $0<c_1<1<c_2<\infty$ and $r_0
> 0$  such that
\begin{equation}\label{gmu}
c_1\ r^{\gamma} \ \le \  \cH^{\gamma}(B\cap K )\, := \,  \mu(B) \
\le \ c_2\ r^{\gamma} \ ,
\end{equation}
for any ball $B=B(x,r)$ with $x\in X$ and $r\le r_0$  -- see for
example \cite{falc,MAT}. Note that  (\ref{gmu}) implies
(\ref{doub}); i.e. the measure $\mu$ is doubling.

 As we shall soon see, the following theorem is an important restatement of
the divergent part of Theorem \ref{t1} in terms of the set
$W^*_{\cal A}(\psi) $ and the measure $\mu$. %with ${\cal H}^f$
%restricted to $\gamma$--dimensional Hausdorff measure ${\cal
%H}^{\gamma}$.

\begin{theorem}\label{t2}
$$ \mu(W^*_{\cal A } (\psi))=  \mu(K)
\qquad {\rm if} \qquad\displaystyle\sum_{n=1}^{\infty}
\left(\psi(3^n) \times \ 3^n \right)^{\gamma}   \ = \ \infty\ \ .
$$
%$$ \cH^{\gamma}(W^*_{\cal A } (\psi) \cap K )=  \cH^{\gamma} (K)
%\qquad {\rm if} \qquad\displaystyle\sum_{n=1}^{\infty}
%\left(\psi(3^n) \times \ 3^n \right)^{\gamma}   \ = \ \infty\ \ .
%$$
\end{theorem}

Note that by definition we have that $$ \mu (W^*_{\cal A } (\psi))
\  = \ \cH^{\gamma}(W^*_{\cal A } (\psi) \cap K ) \qquad {\rm and
\ } \qquad \mu(K) \ = \ \cH^{\gamma}(K) \ .
$$

\noindent The statement of Theorem \ref{t2} is the precise
analogue of the standard Duffin-Schaeffer conjecture for the set
$W^*_{\cal A } (\psi) \cap K $.  The following result enables us
to reduce the proof of the divergent part of Theorem \ref{t1} to
that of establishing Theorem \ref{t2}.

\begin{theorem}\label{t3}
$$ {\it  Theorem  \ \ref{t2}  }  \ \Longrightarrow \ {\it Theorem \
\ref{t1}  \ (divergent \ part) }$$
\end{theorem}

This theorem is a simple consequence of the Mass Transference
Principle.

\vspace*{2ex}

\noindent{\em Proof of Theorem \ref{t3}. \ } Without loss of
generality assume that $\psi(3^n) \to 0 $ as $n \to \infty $.
Otherwise, $W^*_{\cal A } (\psi) = [0,1] $ and the statement is
obvious. We are given that $r^{-\gamma}f(r)$ is monotonic and that
$$
\displaystyle\sum_{n=1}^{\infty} \ f\!\left(\psi(3^n)\right)  \times
\ (3^n)^{\gamma}   \ = \ \infty\ \  .
$$
Let $ \theta : r \to \theta(r):= f(\psi(r))^{\frac{1}{\gamma}}$.
Then,
$$
\displaystyle\sum_{n=1}^{\infty} \left(\theta(3^n) \times \ 3^n
\right)^{\gamma}   \ = \ \infty\ \ .
$$
Thus, Theorem \ref{t2} implies that $\cH^{\gamma}\big(B \cap
(W^*_{\cal A } (\theta) \cap K ) \big)=  \cH^{\gamma} (B) $ for any
ball $B$ in $K$. It now follows via the Mass Transference Principle
that $\cH^{f}\big(B \cap (W^*_{\cal A } (\psi) \cap K ) \big)=
\cH^{f} (B) $ for any ball $B$ in $K$. In particular, this implies
that $\cH^{f}\big(W^*_{\cal A } (\psi) \cap K )= \cH^{f} (K) $ which
together with (\ref{imp}) completes the proof of Theorem \ref{t3}.

\vspace*{-4ex} \hfill $\spadesuit$

\vspace*{2ex}

The upshot  is that  divergent part of Theorem \ref{t1} is a
consequence of Theorem \ref{t2}.

 \vspace*{2ex}

\noindent {\em Remark. \,  } For any ball $B$ in the unit interval
$[0,1]$ and $n \in N$, we have that
\begin{equation}
\label{ubqsep} \mu \Big( \, B  \ \cap \  \bigcup_{0 \leq p \leq
3^n } B( \textstyle{\frac{p}{3^n}},  %\textstyle{\frac{1}{2} }
\ 3^{-n} ) \Big) \ = \ \mu(B) \  \qquad (\mu := {\cal
H}^{\gamma}|_{\substack{ ~ \\ K }})  \ \ .
\end{equation}
This simply makes use of the fact that the distance between
consecutive rationals with fixed denominator $q$ is $1/q$. In view
of (\ref{ubqsep}), if we assume  that $\psi$ is monotonic then the
divergent part of Theorem \ref{t1} is easily seen to be a
straightforward consequence of the local $m$-ubiquity results
established in \cite{BDV03}. Note that Theorem \ref{t1} under the
assumption that $\psi$ is monotonic is enough to determine
Corollaries \ref{c1} and \ref{c2}. Recall, that Mahler's assertion
trivially follows from the dimension part of Corollary \ref{c1}.
In fact, (\ref{ubqsep}) together with the  $m$-ubiquity result
established in \cite{gang4} some   ten years ago  is already
enough to yield the dimension part of Corollary \ref{c1}.

\subsection{Proof of Theorem \ref{t2} \label{prooft2} }

We are given that
\begin{equation}
\label{div2} \displaystyle\sum_{n=1}^{\infty} \left(\psi(3^n)
\times \ 3^n \right)^{\gamma}   \ = \ \infty\ \ .
\end{equation}

\vspace*{2ex}

\noindent {\em Step 1. \,  } We notice that there is no loss of
generality in assuming that $$ \psi(3^n) \leq c \, 3^{-n} \quad {\rm
\ for \ all \ }   n \in \N  {\rm \  and \  }  c
> 0  \  .   $$ Suppose for the moment that this was not case and define
$$\Psi : r \to \Psi(r) := \min \{c/r, \psi(r)\} \ . $$ In view of
(\ref{div2}), it is easily verified that
$$ \textstyle{\sum_{n=1}^{\infty} } \left(\Psi(3^n)
\times \ 3^n \right)^{\gamma}   \ = \ \infty\ \ . $$
 Furthermore,
$ W^*_{\cal A } (\Psi) \subset  W^*_{\cal A } (\psi)$ and so it
suffices to establish  Theorem \ref{t2} for $\Psi$. In particular,
there is no loss of generality in assuming that
\begin{equation}
\label{small} \psi(3^n) \ < \  \textstyle{\frac{1}{2}} \ 3^{-n}
\qquad {\rm for \ all \ } \quad n \in \N  \ \ .
\end{equation}

\vspace*{2ex}

\noindent {\em Step 2. \,  }Let $B$ be an arbitrary ball centered at
a point in $K$ such that $\mu(2B) $ satisfies (\ref{gmu}).
Trivially, $\mu(B) \asymp \mu(2B)$.  The aim is to show that
\begin{equation}
\label{aim}
 \mu(W^*_{\cal A } (\psi) \cap B ) \ \geq \   \mu(B)/C
\ ,
\end{equation}
where $C>0$ is a constant independent of $B$. Theorem \ref{t2} is
then a consequence of Lemma \ref{lem1a}.

 Let $r(B)$ denote the
radius of $B$. Throughout,  $t_0:=t_0(B) $ is a sufficiently large
integer so that
\begin{equation}
\label{t0}
 3^{-t_0} \ < \ r(B) \ \ .
\end{equation}

\noindent For $n \in N$, let
\begin{equation*}
A^*_n(B)   \ :=   \ \bigcup_{\substack{0 \ \leq \ p \  \leq \ 3^n
: \\ ~ \\  (p,3^n) =1  } }^{\circ} \!\!\! B(
\textstyle{\frac{p}{3^n}}, \ \psi(3^{n}) ) \ \cap \ B
\end{equation*}

\noindent The fact that the above union is disjoint is a
consequence of (\ref{small}).  Furthermore, let $B^*_n(r)$ denote
a generic ball centered at a reduced rational $p/3^n$ $ (0  \leq p
\leq 3^n)$ and radius $2r \leq  3^{-n} \! $.  By considering the
$n$'th level $L_n$ of the Cantor set construction in which there
are $2^n$ intervals $I_n$ of common length $3^{-n}\!$, it is
easily verified that if $B^*_n(r) \cap  K \neq \emptyset $  then
its center $p/3^n \in K$; i.e. the rational $ p/3^n $ is an end
point of some interval $I_n$ in $L_n$. Now the measure $\mu$ is
supported on $K$ and satisfies (\ref{gmu}). Hence, for $n > t_0$
$$
\mu \big(B^*_n(r)  \big)  \asymp r^{\gamma}  \qquad {\rm if \ }
\qquad
 B^*_n(r)  \cap  K  \neq \emptyset \ .
$$
Also, it is easily verified that for $n > t_0$
\begin{eqnarray}
\# \{B^*_n(r) \subset B: B^*_n(r)  \cap  K  \neq \emptyset \} &
\asymp & \# \{B^*_n(\textstyle{\frac12} 3^{-n}) \subset B: B^*_n(
\textstyle{\frac12} 3^{-n}) \cap K \neq \emptyset \}  \nonumber \\ & &  \nonumber \\
& \asymp & \frac{ \mu(B)}{ (3^{-n})^{\gamma}} \ . \label{L}
\end{eqnarray}
%Note that in view of (\ref{small}) and the fact that the measure
%$\mu $ is supported on $K$, we have that  $$ \mu \big(A^*_n(B)
%\big) \ = \ \mu \Big(\bigcup_{\substack{0 \ \leq \ p \  \leq \ 3^n
%: \\ ~ \\  (p,3^n) =1  } }^* \!\!\! B( \textstyle{\frac{p}{3^n}},
%\ \psi(3^{n}) ) \ \cap \ B \Big)  \ \ $$ where the $*$ above the
%union indicates the restriction that the rational  $p/3^n$ is in
%$K$.
In view of the above discussion regarding $B^*_n(r)$, a straight
forward geometric argument yields that for $n > t_0$

\begin{equation}
\label{muA}
 \mu(A^*_n(B))    \ \asymp    \  \frac{ \mu(B)}{ (3^{-n})^{\gamma}} \times
 \left(\psi(3^n)
 \right)^{\gamma}  \  = \  \mu(B) \; \left(\psi(3^n)
\times \ 3^n \right)^{\gamma}  \ \ .
\end{equation}

\noindent This together with (\ref{div2}) implies that
\begin{equation}
\label{muAdiv}
 \sum_{n=1}^{\infty} \mu(A^*_n(B))   = \infty  \ .
\end{equation}
%\noindent
Finally, note that

$$
W^*_{\cal A } (\psi) \cap B = \limsup_{n \to \infty} A^*_n(B) \ .
$$

\vspace*{2ex}

\noindent {\em Step 3. \,  }The key to establishing Theorem \ref{t2}
is the following `local' pairwise quasi--independence result.

\begin{lemma}[Local pairwise quasi--independence] \label{In}
%There exist constants $t_0:=t_0(B)$ and $C >0$  such that for all
%$\, n
%> m > t_0 $,
There exists a constant $C >0$  such that for all $\, n
> m > t_0 $,
\begin{equation}
\label{quasienq}
 \mu (A^*_m(B) \cap A^*_n(B)) \  \leq  \
\frac{C}{\mu(B)} \;  \mu (A^*_m(B)) \  \mu( A^*_n(B))  \ \  .
\end{equation}
\end{lemma}

Hence, (\ref{muAdiv}) and Lemma \ref{In} together with Lemma
\ref{lem2a} implies (\ref{aim}). This completes the proof of
Theorem \ref{t2} assuming of course the local pairwise
quasi--independence result -- this we now prove.

\subsection{Proof of Lemma \ref{In}: Local pairwise
quasi--independence}

Recall, that  $B$ is some fixed  ball centered at a point in $K$
such that $\mu(2B) $ satisfies (\ref{gmu}) and $t_0:=t_0(B) $ is
chosen sufficiently large so that (\ref{t0}) is satisfied.
% Let $r(B)$ denote the
%radius of $B$. Choose $t_0:=t_0(B) $ sufficiently large so that
%$$
%3^{-t_0} \ < \ r(B) \ \ . $$
%Then, for any $m > t_0$ there exist at
%least two intervals  $I_m$ at the $m$'th level $L_m$ of the Cantor
%set construction which are contained within $B$.

Fix a pair $n$ and $m$ with $n > m > t_0$. We proceed by  consider
two cases depending on the size of $\psi(3^n)$ compared to $
\textstyle{\frac{1}{2}} \, 3^{-n} $.

\vspace*{2ex}

\noindent{\em Case (i):  $n > m > t_0 $ such that  $ 3^{-n}
 \geq  2 \, \psi(3^m)$. \,  } Fix some ball $B(
p/3^m \! , \psi(3^{m}) ) $ where $0 \leq p \leq 3^m$. It is easily
verified on assuming (\ref{small}), that %under the condition imposed
%in  case (i)
$$
B \big( \textstyle{\frac{p}{3^m}}, \psi(3^{m}) \big) \ \cap \ B\big(
\textstyle{\frac{t}{3^n}}, \psi(3^{n}) \big) \ = \ \emptyset  \
\quad {\rm \ for \ any \ } \ 0 \leq t \leq 3^n \ {\rm \  with \ } \
(t, 3^n) =1 \ .
$$
%for any $0 \leq t \leq 3^n$ with $(t, 3^n) =1 $.
In view of this,
$$
A^*_m(B) \cap A^*_n(B) \  = \ \emptyset \ . $$ Thus, $ \mu (A^*_m(B)
\cap A^*_n(B))  =  0 $ and (\ref{quasienq}) is trivially satisfied
for any constant $C \geq 0 $.

\vspace*{2ex}

\noindent{\em Case (ii):  $n > m > t_0$ such that  $
 3^{-n}
 < 2 \, \psi(3^m)$. \,  } For the sake of clarity,  let us write  $B^*_n(\psi)$
 for the generic ball $B^*_n(\psi(3^n))$  and simply  $ B^*_n$ for
 the generic ball $B^*_n(\textstyle{\frac12} \, 3^{-n})$. Recall,
 that by definition a generic ball $B^*_n(r)$ is  centered at a reduced  rational
  $p/3^n$ $ (0 \leq p \leq 3^n)$ and if it  has non-empty
 intersection with $K$ then   $p/3^n \in K $ -- see Step 2 in
 \S\ref{prooft2}.

It is easily verified  that

\begin{eqnarray}
\mu (A^*_m(B) \cap A^*_n(B))  & := & \mu \Big( \big(
\bigcup_{\substack{0 \ \leq \ p \  \leq \ 3^m : \\ ~ \\  (p,3^m)
=1 } }^{\circ} \!\!\! B( \textstyle{\frac{p}{3^m}}, \ \psi(3^{m})
) \ \cap \ B \big) \ \cap \ A^*_n(B) \Big) \nonumber \\ & &
\nonumber
\\ & \leq & \sum_{\substack{0 \ \leq \ p \  \leq \ 3^m \, : \,
(p,3^m) =1,
\\  ~  \\  B( \frac{p}{3^m},  \psi(3^{m}) ) \, \cap \,
B \,  \cap \,   K \neq \emptyset} } \mu \big( B(
\textstyle{\frac{p}{3^m}}, \ \psi(3^{m}) ) \ \cap \  A^*_n(B)
\big) \nonumber
\\ & &
\nonumber \\ & \leq & {\cal N}(m,B) \ \times \ \mu \big( B^*_m(
 \psi ) \, \cap \,  A^*_n(B)
\big) \ , \label{q1}
\end{eqnarray}
where
$$
{\cal N}(m,B) \ := \ \# \{ B^*_m( \psi ) \, : \, B^*_m( \psi ) \,
\cap \,  B \, \cap \,  K  \, \neq \, \emptyset \} \ .
$$

\noindent In view of  the fact that $m > t_0 $ and $t_0$ satisfies
(\ref{t0}), we have that
\begin{eqnarray}
{\cal N}(m,B) &  \stackrel{(\ref{small})}{ \leq  }  &  \# \{
B^*_m( \psi ) \subset 2B  \, : \,
B^*_m( \psi ) \, \cap \,  K  \, \neq \, \emptyset \} \nonumber \\
& & \nonumber \\
& \stackrel{(\ref{L})}{ \leq  } & c_3 \ \mu(B) \ \times \
(3^m)^{\gamma} \ \ , \label{q2}
\end{eqnarray}
where $c_3 > 0 $ is a constant. We now obtain an upper bound for
$\mu \big( B^*_m( \psi ) \, \cap \,  A^*_n(B) \big)$. Without loss
of generality, we assume that $ B^*_m( \psi ) \, \cap \,  B \,
\cap \,  K  \, \neq \, \emptyset $ since otherwise $\mu \big(
B^*_m( \psi ) \, \cap \,  A^*_n(B) \big) =0 $ and there is nothing
to prove. For a fixed generic ball $B^*_m( \psi ) $, a relatively
simple geometric argument yields that
\begin{eqnarray*}  \# \{
B^*_n( \psi )  \, : \, B^*_n( \psi ) \, \cap \, B^*_m( \psi ) \, K
\, \neq \, \emptyset \}  & \leq & \# \{ B^*_n  \, : \,
B^*_n \, \cap \, B^*_m( \psi ) \, \cap \, K  \, \neq \, \emptyset \} \nonumber \\
& & \nonumber \\
& \stackrel{(\ref{gmu})}{ \leq  } & {\textstyle{\frac{c_2}{c_1} }}
\ \big( \psi(3^m) \,  \, 3^{-n} \big)^{\gamma} \  + \ 2 \ .
\nonumber
\end{eqnarray*}
The  `plus $2$' term above simply accounts for `edge effects'.
Hence,
\begin{eqnarray}
\mu \big( B^*_m( \psi ) \, \cap \,  A^*_n(B) \big)  & \leq   & \mu
\big(B^*_n( \psi ) \big) \, \times \, \big( \,
{\textstyle{\frac{c_2}{c_1} }} \ \big( \psi(3^m) \,  \, 3^{-n}
\big)^{\gamma} \  + \ 2 \, \big) \nonumber \\ & & \nonumber \\  &
\stackrel{(\ref{gmu})}{ \leq  } & {\textstyle{\frac{c_2^2}{c_1} }}
\ \big( \psi(3^n) \ \psi(3^m) \, \, 3^{-n} \big)^{\gamma} \  + \ 2
\, c_2 \big( \psi(3^n) \big)^{\gamma} \  \ .
 \label{q3}
\end{eqnarray}
On combining (\ref{q1}), (\ref{q2}) and (\ref{q3}), we obtain that
\begin{eqnarray*}
\mu \big( A^*_m(B) \, \cap \,  A^*_n(B) \big)  & \leq   &  c_2 \,
c_3 \ \mu(B) \, \,   \big( \, 3^n \, \psi(3^n) \ 3^m \, \psi(3^m)
\, \big)^{\gamma} \, \times \, \big( {\textstyle{\frac{c_2}{c_1}
}} \,  + \, 2 \,  \big( 3^n \psi(3^m)
\big)^{- \gamma} \big) \\ & & \\
& \stackrel{case \,(ii)}{ \leq  } &  c_4 \ \mu(B) \, \, \big( \,
3^n \, \psi(3^n) \ 3^m \, \psi(3^m) \, \big)^{\gamma}
\\ & & \\
& \stackrel{(\ref{muA})}{ \leq  } &   \frac{C}{\mu(B)} \, \mu
\big( A^*_m(B) \big) \  \mu \big( A^*_n(B) \big) \ ,
\end{eqnarray*}
where $c_4>0$ and $C>0$ are absolute constants. Thus,
(\ref{quasienq}) is satisfied.

\vspace*{2ex}

On combining the above two cases concludes the proof of the local
pairwise quasi--independence statement -- Lemma \ref{In}.

\vspace*{-4ex} \hfill $\spadesuit$

\section{General `missing digit' sets $K_{J(b)}$}

Let $b \geq 3 $ be an integer and let $J(b)$ be a proper subset of
$ S := \{0,1, \ldots, b-1 \}$ with $\# J(b) \geq 2 $. Furthermore,
let $K_{J(b)}$ denote the set of real numbers in the unit interval
$[0,1]$ whose base $b$ expansions consist exclusively of digits
within $J(b)$. Equivalently, $K_{J(b)}$ is the set of  real
numbers in the unit interval whose base $b$ expansions are free of
the digits in $S \setminus J(b)$. Thus, $K_{J(b)}$ is a natural
generalization of the middle third Cantor set $K$ -- simply put
$b=3$ and let $J(b) = \{0,2\}$.
%\vspace*{2ex}
Naturally, one can ask whether Mahler's assertion remains valid
with the middle third Cantor set $K$ replaced by the general
Cantor set $K_{J(b)}$.

\vspace*{1ex}

It is easily versified that
$$
 \dim K_{J(b)} \ =  \gamma^*  \ := \ \frac{\log \# J(b) }{\log b } \
 .
$$
Moreover, there exists a finite measure $\mu $ supported on
$K_{J(b)}$ such that $$ {\cal H}^{\gamma^*}  (B \cap K_{J(b)} ) \
:= \ \mu (B) \ \asymp \ r^{\gamma^*} \ \ ,  $$ for any ball
$B=B(x.r)$  with $x \in K_{J(b)}$ and $r \leq r_0 $. Both the
dimension and measure statements above can be deduced from
standard results in fractal geometry -- see for example
\cite{falc,MAT}.

Now let  ${\cal A}(b):= \{ b^n : n = 0,1,2, \ldots \} $ and
consider the set
$$
W_{{\cal A}(b)}(\psi):= \{x \in [0,1] :|x - p/q| < \psi(q) \ {\rm
for\ i.m.\ } \; (p,q) \in \Z \times {\cal A}(b) \} \ . $$ The
arguments involved in establishing Theorem \ref{t1}  can be
modified in the obvious manner to yield the following
generalization of Theorem \ref{t1}.

\begin{theorem}\label{t11}
Let $f$ be a dimension function such that $r^{-\gamma^*}f(r)$ is
monotonic.  Then
$$
\cH^f(W_{{\cal A }(b)} (\psi) \cap K_{J(b)} )=\left\{
\begin{array}{cl}
0& {\rm if}  \qquad\displaystyle \sum_{n=1}^{\infty} \
f\!\left(\psi(b^n)\right)  \times \ (b^n)^{\gamma^*} \ < \ \infty\, \\[4ex]
\cH^f(K_{J(b)})  & {\rm if} \qquad\displaystyle\sum_{n=1}^{\infty} \
f\!\left(\psi(b^n)\right) \times \ (b^n)^{\gamma^*}   \ = \ \infty\,
\end{array}
\right..
$$
\end{theorem}

Simple consequences of the theorem are the following corollaries
which clearly imply the existence of very well approximable numbers,
other than  Liouville numbers, in the  Cantor set $K_{J(b)}$. In
fact, it follows  that  $$ \dim ( \, ({\cal W}\setminus {\cal L
})\cap K_{J(b)} )  \ \geq \ \gamma^*/ 2 \  . $$
%\ \geq \ \gamma^*/ 2 $
%i.e.
%$$ ({\cal W}\setminus {\cal L })\cap K_{J(b)} \ \neq  \ \emptyset \ .
%$$
%In fact, $ \dim ( \, ({\cal W}\setminus {\cal L })\cap K_{J(b)} )
%\ \geq \ \gamma^*/ 2 $ and so $({\cal W}\setminus {\cal L })\cap
%K_{J(b)} \ \neq  \ \emptyset$.
\begin{corollary}\label{c11} For $\tau \geq 1 $,
$ \dim (W_{{\cal A }(b)} (\tau) \cap K_{J(b)} ) \ = \ \gamma^*/
\tau  $.  In particular, for $\tau > 1 $
$$\cH^{\gamma^*/\tau} (W_{{\cal A }(b)} (\tau) \cap K_{J(b)} ) \ = \ \infty \ . $$
\end{corollary}

\begin{corollary}\label{c21} For $\alpha \geq 2 $,
$ \dim ( E(\alpha) \cap K_{J(b)} ) \ \geq \ \gamma^*/ \alpha  $.
In particular, for $\alpha > 2 $
$$\cH^{\gamma^*/\alpha} ( E(\alpha) \cap K_{J(b)}) \ = \ \infty \ . $$
\end{corollary}

\section{Concluding Remarks \label{conc}}

\subsection{ Mahler's assertion -- explicit examples \label{explicit}}
For any real number  $\tau > 2 $, consider the irrational number
$$
\xi \ = \  \xi(\tau) \ := \  2 \, \sum_{n=1}^{\infty}   \,
3^{-\tau_n} \ \qquad {\rm where } \ \qquad  \tau_n :=  [\tau^n] \
.
$$
Clearly $\xi$ is irrational since its base 3 expansion is not
periodic.  Clearly $\xi $ lies in $K$ since only the digits $0$
and $2$ appear in its base $3$ expansion. The following result
implies that
$$\xi \ \in \  ({\cal W}\setminus {\cal L })\cap K \ .
$$

\noindent Thus, for each $\tau > 2$  the  explicit irrational
number $\xi= \xi(\tau)$ satisfies  Mahler's assertion.

\begin{lemma} (i) If $\, \tau \geq \frac{1}{2} (\sqrt{5} +3) $, then
$$
\xi \ \in \  E(\tau)  \ \cap \ K \ \ .
$$
(ii) If $\, 2< \tau <  \frac{1}{2} (\sqrt{5} +3)$, then for any
$\epsilon
> 0$
$$
\xi \ \in \ \big( \, W(\tau-\epsilon) \setminus W\big( \textstyle{
\frac{2 \tau -1}{\tau -1} }  + \epsilon \big) \, \big) \ \cap \ K
\ \ .
$$
\label{exp}
\end{lemma}

\vspace*{2ex}

\noindent{\em Remark.  } Note that we are only able to conclude
the stronger  exact order statement (part (i)) under the
assumption that $\tau \geq \frac{1}{2} (\sqrt{5} +3) $. However,
statements of this type are reminiscent of numerous results in
transcendence theory; see, for example \cite[Theorems 7.7 \&
8.8]{bug}. By the definition of the exact order set $E(\tau)$, we
trivially have that  $\xi \in {\cal W}\setminus {\cal L }$.
Regarding part (ii) of the above lemma, by choosing $\epsilon > 0
$ sufficiently small so that $ \tau - \epsilon > 2 $,  we also
have that $\xi \in {\cal W}\setminus {\cal L }$. %Finally notice
%that  the lemma yields uncountably many explicit irrational
%numbers satisfying Mahler's assertion.

\vspace*{1ex}

\noindent{\em Proof of Lemma \ref{exp}. \ } As already mentioned
above, the fact that $\xi \in K$ is trivial.  For $s \in \N $, let
$$
q_s \ :=  \  3^{{\tau}_s}   \quad { \rm and \ } \quad  p_s  \ := \
 q_s  \, \times \  2 \, {\displaystyle{\sum_{n=1}^{s}}}  \, 3^{-{\tau}_{n}} \ .
$$
It is easily verified that for all $s \in \N$, $(p_s,q_s) = 1 $
and that
\begin{equation}
\label{qs}
 \textstyle{\frac13} \, {q_s}^{\! \tau}  \ < \ q_{s+1}  \ < \  3^{\tau} \,  {q_s}^{\!
\tau}  \  \ .
\end{equation}
Furthermore,
\begin{equation}
\label{eeg}
 \frac{2}{3^\tau} \; \frac{1}{{q_s}^{\! \tau}}  \ \stackrel{(\ref{qs})}{<} \
 \frac{2}{q_{s+1}} \ < \ \Big\vert \xi- \frac{p_s}{q_s}
\Big\vert  \ <  \ \frac{3}{q_{s+1}} \ \stackrel{(\ref{qs})}{<} \
\frac{9}{{q_s}^{\! \tau}} \  \ .
\end{equation}

\noindent Fix some  $\epsilon > 0$. Since $9/{q_s}^{\!\! \tau} <
1/{q_s}^{\!\! \tau - \epsilon} $ for all sufficiently large $s$,
we have  that $$ \xi \in W(\tau-\epsilon)  \ .  $$ By assumption,
$ \tau > 2 $. Therefore, there exists an integer $s_{\rm o}$, such
that for all $ s \geq s_{\rm o} $
$$
\Big\vert \xi- \frac{p_s}{q_s} \Big\vert  \
\stackrel{(\ref{eeg})}{ < }   \ \frac{9}{{q_s}^{\! \tau}} \ < \
\frac{1}{{2 \, q_s}^{\! 2}} \ \ .
$$
It follows, via a standard result in the theory of continued
fractions (Legendre's theorem), that for each  $ s \geq s_{\rm o}$
the rational $p_s/q_s $ is a convergent to $\xi$. Now, with a
slight abuse of notation, it is well known that  succussive
convergents $p_n/q_n $ and $ p_{n+1} /q_{n +1} $ to any irrational
number $x$ lie on either side of $x$ and that
\begin{equation}
\label{cf}
 \frac{1}{q_n(q_n +q_{n+1})} \ < \ \Big\vert x - \frac{p_n}{q_n}
\Big\vert  \ <  \ \frac{1}{q_n q_{n+1}} \   \ .
\end{equation}
Also, the denominators $\{ q_n :  n \geq 2 \} $ of the convergents
form a strictly increasing sequence. For these standard statements
from the theory of continued fractions the reader is referred to
\cite[\S1.2]{bug}. By construction or rather by definition, we
have that
$$
 \frac{p_s}{q_s} \ < \ \frac{p_{s+1} }{q_{s+1}}  \ < \
\xi   \ \ . $$ Thus, for $ s \geq s_{\rm o} $ the convergents
$p_s/q_s$ and $ p_{s+1}/q_{s+1} $ are not succussive convergents
to $\xi$ and so there exists at least one other convergent with
denominator between $q_s$ and $ q_{s+1}$.  With this in mind, let
$p_*/q_*$ be the next convergent to $\xi$ after $p_s/q_s $. Thus,
$\frac{p_*}{q_*} \neq \frac{p_{s+1} }{q_{s+1}}  $  and $q_s < q_*
< q_{s+1} $.  In view of (\ref{cf}), we have that
\begin{equation*}
\label{cf1}
 \frac{1}{q_s(q_s +q_{*})} \ < \ \Big\vert \xi- \frac{p_s}{q_s}
\Big\vert  \ <  \ \frac{1}{q_s q_{*}} \   \ .
\end{equation*}
This  together with (\ref{eeg}) and the fact that $\tau > 2$,
implies that there exists an integer $ s_{ \rm 1} \geq s_{\rm o}$
such that  for $s \geq s_{ \rm 1}$
\begin{equation}
\label{q*} \textstyle{\frac{1}{10} } \, q_s^{\tau-1}  \ < \ q_* \
< \ \textstyle{\frac{\ 3^{\tau}}{2} } \; q_s^{\tau-1} \   \ .
\end{equation}

Now suppose there exists a rational $p/q $ and $s \geq  s_{\rm 1}$
such that $q_s <  q < q_{s+1}$ and
\begin{equation}
\label{pq}
  \Big\vert \xi- \frac{p}{q}
\Big\vert  \ <  \ \frac{1}{q^{\tau + \epsilon}} \
% \qquad {\rm and } \qquad  q_s <  q < q_{s+1}   \
\ .
\end{equation}
Thus, $p/q $ is a convergent to $\xi$ with
\begin{equation}
\label{q>}
 q \ \geq \ q_*  \ \ .
\end{equation}

Now let $\tilde{p}/\tilde{q}$ be the next convergent to $\xi$
after $p/q$. Then,
\begin{equation}
\label{qtil<}
 \tilde{q}  \ \leq  \ q_{s+1}  \ \ .
\end{equation}
 In view of (\ref{cf}), we have that
\begin{equation*}
\label{cf2}
 \frac{1}{q(q + \tilde{q})} \ < \ \Big\vert \xi- \frac{p}{q} \Big\vert  \ <
\ \frac{1}{q \, \tilde{q}} \   \ .
\end{equation*}
This together with (\ref{pq}) implies that
\begin{eqnarray}
\label{qtil}
  \tilde{q}   \   >  \ \textstyle{ \frac{1}{2} } \,  q ^{\tau + \epsilon -1 }
  \ & \stackrel{(\ref{q>})}{\geq} & \
  \textstyle{ \frac{1}{2} }   \, q_* ^{\tau + \epsilon -1 }  \nonumber
  \\& & \nonumber \\
 \ & \stackrel{(\ref{q*})}{>} & \  \textstyle{ \frac{1}{2} } \,
 \big(\textstyle{ \frac{1}{10} }\big)^{\tau + \epsilon -1 }
   \ \big(q_s^{\tau-1} \big)^{\tau + \epsilon -1 }   \ .
\end{eqnarray}

Suppose for the moment that $\tau \geq \frac{1}{2} (\sqrt{5} +3)$.
Then,
\begin{equation}
\label{te}
 (\tau-1)(\tau + \epsilon -1 ) \ > \ \tau \ \ .
\end{equation}  It follows that  there
exists an integer $ s_{ \rm 2} \geq s_{\rm 1}$ such that for $ s
\geq s_{ \rm 2}$
$$
\textstyle{ \frac{1}{2} } \, \big(\textstyle{ \frac{1}{10}
}\big)^{\tau + \epsilon -1 } \ \big(q_s^{\tau-1} \big)^{\tau +
\epsilon -1 } \ > \ 3^{\tau} q_s^{\tau}  \
\stackrel{(\ref{qs})}{>} \ q_{s+1} \ \ .
$$
This together with (\ref{qtil}) implies that  $\tilde{q} >
q_{s+1}$ which in view of  (\ref{qtil<}) is a contradiction.  The
upshot is that there are at most finitely many rationals $p/q$
satisfying  inequality (\ref{pq}).  Hence,
$$ \xi \notin  W(\tau+\epsilon)  \  $$
and this completes the proof of part (i) of the lemma. Part (ii)
of the lemma follows from the observation that if
$$
\epsilon \ > \ \frac{\tau}{\tau -1 }  - \tau + 1  \ > \ 0  \ \ ,
$$
then (\ref{te}) is  satisfied  and we are still able to force the
contradiction that $\tilde{q} > q_{s+1}$.

\vspace*{-3ex} \hfill $\spadesuit$

%%%%%%%%%%%%%%%%%%%%%%%%%%%%%%%%%%%%%%%%%%%%%%%%%%%%%%%%%%%%%%%%%%%%%%%%%%%%%%%%%%%%%%%%%%%
%%%%%%%%%%%%%%%%%%%%%%%%%%%%%%%%%%%%%%%%%%%%%%%%%%%%%%%%%%%%%%%%%%%%%%%%%%%%%%%%%%%%%%%%%%%
%%%%%%%%%%%%%%%%%%%%%%%%%%%%%%%%%%%%%%%%%%%%%%%%%%%%%%%%%%%%%%%%%%%%%%%%%%%%%%%%%%%%%%%%%%%
%%%%%%%%%%%%%%%%%%%%%%%%%%%%%%%%%%%%%%%%%%%%%%%%%%%%%%%%%%%%%%%%%%%%%%%%%%%%%%%%%%%%%%%%%%%
%%%%%%%%%%%%%%%%%%%%%%%%%%%%%%%%%%%%%%%%%%%%%%%%%%%%%%%%%%%%%%%%%%%%%%%%%%%%%%%%%%%%%%%%%%%

\vspace*{1ex}

\noindent{\em Remark.  } Fix some $ \tau \geq \frac{1}{2}
(\sqrt{5} +3)$. Let $\lambda > 0$ be a real number.  On adapting
the above argument in the obvious manner, it is readily verified
that
$$
\xi = \xi( \tau, \lambda)  \ := \ \  2 \, \sum_{n=1}^{\infty}   \,
3^{- [\lambda \tau^n]} \   \ \ \in \ \  E(\tau) \cap K  \ \ .
$$
This clearly yields uncountably many explicit irrational numbers
in $E(\tau) \cap K$.

\subsection{ Mahler's assertion -- what do we really expect? \label{expect}} We
suspect that  lower bound estimate (\ref{walk})  for $\dim ( \,
({\cal W}\setminus {\cal L })\cap K  ) $ is far from the truth. It
is highly likely that:
\begin{equation} \label{walk2}
\dim ( \, ({\cal W}\setminus {\cal L })\cap K  )  \ = \ \gamma \
:= \ \dim K \ \  .
\end{equation}

Recall, estimate (\ref{walk}) is obtained by considering rationals
with denominators restricted to the set ${\cal A}:= \{ 3^n : n =
0,1,2, \ldots \} $ and showing that  $ \dim (W_{\cal A } (\tau) \cap
K )  = \gamma/ \tau $  for $\tau \geq 1$ -- Corollary \ref{c1}.  The
restriction of the denominators to the set ${\cal A} $ is absolutely
paramount to most of the arguments employed in this paper -- it
forces the rationals of interest to lie in $K$. For a brief moment,
let us forget about intersecting well approximable sets with $K$. It
is well known that $ \dim W_{\cal A } (\tau) = 1/\tau $ $ (\tau \geq
1)$ where as $ \dim W(\tau) = 2/\tau $ $(\tau \geq 2)$ -- see
\cite[\S12.5]{BDV03} and \cite[Chp.10]{Harman}. Thus, by restricting
the denominators to the set ${\cal A}$, the dimension is reduced or
rather re-scaled by a factor of $1/2$. It is therefore reasonable to
speculate that this same scaling is present when considering the
dimensions of the sets $W_{\cal A } (\tau) \cap K$ and $W(\tau) \cap
K$. This leads us to the following statement which would imply
(\ref{walk2}) in the same way that Corollary \ref{c1} is used to
establish (\ref{walk}).

\vspace*{0.4ex}

{\bf Statement 1. \, } For $\tau \geq 2 $, $\ \ \dim (W(\tau) \cap
K) \, = \, 2 \, \gamma/ \tau $.

\vspace*{0.4ex}

\noindent{\em Remark. }It is reasonably easy to obtain the upper
bound estimate; namely $\dim (W(\tau) \cap K) \leq 2 \, \gamma/
\tau $. It also follows directly from Corollary 2 of \cite{pv}.
Thus the  problem lies in establishing the complementary lower
bound estimate.

  Even if the above statement turns out to be false, we
have  every reason to believe the following weaker statement which
would still imply (\ref{walk2}).

\vspace*{0.4ex}

{\bf Statement 2. \, } Let  $\tau > 2 $. Then, $ \ \dim (W(\tau)
\cap K) \, \to \, \gamma  \, $ as $ \,   \tau \to 2  $.

\vspace*{0.4ex}

\noindent Ideally, one would like to obtain a complete metric
theory for the set $ W(\psi) \cap K$. The above statements  would
then  be simple corollaries  of such a theory. Equivalently, one
would like to obtain the analogue of Theorem \ref{t1} for the set
$ W(\psi) \cap K$.
 Without, imposing the condition
that $\psi$ is monotonic the problem seems harder than
establishing the Duffin-Schaeffer conjecture.  It is quite likely
that the necessary ideas and techniques required  in proving the
analogue of Theorem \ref{t1} for the set $ W(\psi) \cap K$ would
also lead to a proof of the Duffin-Schaeffer conjecture.   In view
of this, suppose we impose the condition that $\psi$ is monotonic.
It is then  relatively straightforward  to obtain the following
convergent result:
$$
\cH^f(W(\psi) \cap K ) \ = \  0 \qquad {\rm if}
\qquad\displaystyle \sum_{n=1}^{\infty} \ f\!\left(\psi(n)\right)
\times \ n^{2\gamma -1 } \ < \ \infty\ .
$$
When $f : r \to r^s$, the above statement follows directly from
Theorem 2 of \cite{pv}. The ideas in \cite{pv} can be easily
modified to deal with general dimension functions. As far as we
are aware, there has been absolutely no progress towards
establishing a divergent theory even in the case that $\cH^f$ is
the Cantor measure $ \mu :=  \cH^{\gamma}|_K$.

It is worth mentioning that progress has  recently  been made
towards establishing a  convergent theory for  general fractal sets
supporting so called `friendly' measures. The Cantor set $K$
together with the Cantor measure $\mu$ fall into this general theory
-- see \cite{KLW,pv,Weiss}.

%\noindent $ \bullet\ $
\subsection{Another problem of Mahler concerning Cantor sets}
It is believed  by most experts that the base b expansion of  any
irrational algebraic number $x$ is normal. Restricting to the case
$b=3$, this would mean that each of the digits $0,1$ and $2$ occurs
the expected number of times in the base $3$ expansion of $x$. In
particular, for $N$ sufficiently large we expect the digit $1$ to
occur around $N/3$ times in the first $N$ terms of the expansion of
$x$ in base 3. However, it is not even  know that the digit $1$
occurs at least once. Equivalently,

\vspace{2ex}

\noindent {\bf Mahler's Problem \cite{mah}. \ }  {\em Are irrational
elements of Cantor's set necessarily transcendental? Thus does
Cantor's set contain no irrational algebraic elements? }

\vspace{2ex}

\noindent{\em Remark.  }  Roth's theorem on rational approximation
to algebraic numbers states that irrational algebraic numbers are
not very well approximable.  Thus, irrational algebraic numbers
are clearly excluded from Mahler's Assertion of \S\ref{prob}.

To our knowledge, it is not even known  whether or not  quadratic
irrationals avoid the middle third Cantor set $K$. It is easy to
see that the golden ratio $(\sqrt{5} -1)/2$ is not in $K$. Let
$[a_1, a_2, a_3, \ldots ]$ be  the standard continued fraction
expansion of a real number $x$ in [0,1]. If $ a_1=a_2 = 1$, then
it is easily  verified that $x$ lies in the interval $
[{\textstyle{\frac12, \frac23}} )$. However, this interval clearly
misses $K$. Now, $(\sqrt{5} -1)/2 = [\bar{1}]:= [1,1,1, \ldots ] $
and so is not in $K$.

Refining the problem even further, it is not at all obvious (to us
at any rate) that  all `simple' quadratic irrationals  avoid $K$;
i.e. if $x = [\bar{n}]$  then $x \notin K $ for any $n \in \N$. It
may well be the case that Mahler's problem is no easier for
quadratic irrationals. The fact that quadratic irrationals have
periodic continued fraction expansions may  simply be a red
herring!

\vspace*{12mm}

\noindent{\bf Acknowledgments: \, } SV would like to thank the
laughing  girls Mr Tiger Scuba Diver Wiver Tiger (Iona) and
Dorothy Princess Angel Flower Star (Ayesha) for their continuous
gibberish -- all of which is of  the highest quality !! SV would
also like to thank Yann Bugeaud for his detailed
comments/suggestions during the `writing up' period of the paper.
JL would like to thank SV for his persistence and for his near
legendary opening gambit ``Alright me dear? 'ere, I've been
thinking!''

\vspace{5mm}

\noindent Jason Levesley: \ Department of Mathematics, University
of York,

\vspace{-2mm}

 ~ \hspace{19mm}  Heslington, York, YO10 5DD, England.

%\vspace{0mm}

 ~ \hspace{19mm} e-mail: jl107@york.ac.uk

 \vspace{5mm}

\noindent Cem Salp:  \hspace{8mm} Department of Mathematics,
University of York,

\vspace{-2mm}

 ~ \hspace{19mm}  Heslington, York, YO10 5DD, England.

%\vspace{0mm}

 ~ \hspace{19mm} e-mail: cs502@york.ac.uk

\vspace{5mm}

\noindent Sanju L. Velani: Department of Mathematics, University
of York,

\vspace{-2mm}

 ~ \hspace{19mm}  Heslington, York, YO10 5DD, England.

%\vspace{0mm}

 ~ \hspace{19mm} e-mail: slv3@york.ac.uk


\vspace*{5mm}


\begin{thebibliography}{99}

\bibitem{BS} A. Baker and W. M. Schmidt. Diophantine approximation
and Hausdorff dimension. {\em Proc. Lond. Math. Soc.}, (3) 21
(1970) 1-11.


\bibitem{exactBDV} V. Beresnevich, H. Dickinson  and
S. L. Velani : Sets of exact `logarithmic' order in the theory of
Diophantine approxiamtion. {\em Math. Ann.} 321 (2001) 253--273.


\bibitem{BDV03} V. Beresnevich, H. Dickinson and S. L. Velani :
Measure Theoretic Laws for limsup sets, Pre-print (97pp):
arkiv:math.NT/0401118. To appear:  {\em Memoirs of the AMS}.

\bibitem{mtp} V. Beresnevich and S. L. Velani :
A Mass Transference Principle and the Duffin--Schaeffer conjecture
for Hausdorff measures, Pre-print (22pp): arkiv:math.NT/0401118.


\bibitem{bug}
Y. Bugeaud : {\em Approximation by algebraic numbers. } Cambridge
Tracts in Mathematics 160, C.U.P., (2004).
%
%
%\bibitem{det} H. Dickinson and S. L. Velani :
%\emph{Hausdorff measure and linear forms},
% J. reine angew. Math., 490 (1997) 1-36.

\bibitem{gang4}
M. M. Dodson, M. V. Meli$\acute{a}$n, D. Pestana and S. L. Velani
: Patterson measure and Ubiquity. {\em  Ann. Acad. Sci. Fenn.},
20:1 (1995) 37--60.


\bibitem{falc}
K. Falconer : {\em Fractal Geometry : Mathematical Foundations and
Applications.} John Wiley \& Sons, (1990).


\bibitem{Gut}R. G\"{u}ting :
On  Mahler's function $\theta_1$. {\em Michigan Math. J.}, 10
(1963), 161-179.

\bibitem{Harman}
G. Harman : {\em Metric Number Theory.} LMS Monographs 18,
Clarendon Press, Oxford, (1998).




\bibitem{KW}
D. Kleinbock and B. Weiss : Badly approximable vectors on
fractals.  Pre-print (31pp):
http://people.brandeis.edu/~kleinboc/Pub/bad.pdf. To appear: {\em
Israel J. Math.}



\bibitem{KLW}
D. Kleinbock, E. Lindenstrauss and B. Weiss : On fractal measures
and Diophantine approximation.  Pre-print (41pp):
http://people.brandeis.edu/~kleinboc/Pub/friendly.pdf. To appear:
{\em Selecta Mathematica}.

\bibitem{KTV}
S. Kristensen, R. Thorn and S. L. Velani : Diophantine approximation
and  badly approximable sets. Pre-print (40pp):
arkiv:math.NT/0405433. To appear:  {\em Advances in Math}.

\bibitem{mah}
K. Mahler : Some suggestions for further research. {\em Bull.
Austral. Math. Soc.}, 29 (1984), 101-108.

\bibitem{MAT}
P. Mattila : {\em Geometry of sets and measures in Euclidean space},
CUP, Cambridge studies in advance mathematics {\bf 44} (1995)

\bibitem{pv}
A. D. Pollington and  S. L. Velani : Metric Diophantine
approximation and `absolutely  friendly' measures,  Pre-print
(11pp): http://www.arxiv.org/abs/math.NT/0401149. To appear: {\em
Selecta Mathematica}.

\bibitem{Sp}
V. G. Sprind\v{z}uk : {\em Metric theory of Diophantine
approximation} (translated by R. A. Silverman). V. H. Winston \&
Sons, Washington D.C. (1979).


\bibitem{Weiss} B. Weiss : Almost no points on a Cantor set are very
well approximable. {\em Proc. R. Soc. London.} 457 (2001), 949--952.

\end{thebibliography}
\end{document}